\documentclass[12pt]{amsart}
\usepackage{a4,amsthm,epsf}         

\theoremstyle{plain}
\newtheorem{theorem}{Theorem}[section]
\newtheorem{lem}[theorem]{Lemma}
\newtheorem{prop}[theorem]{Proposition}
\newtheorem{cor}[theorem]{Corollary}

\theoremstyle{definition}

\theoremstyle{remark}
\newtheorem*{rem} {Remark}

\def\g{{\mathfrak{g}}}
\def\h{{\mathfrak{h}}}

\def\C{{\mathbb{C}}}
\def\R{{\mathbb{R}}}

\def\nsmallskip{\smallskip\noindent}
\def\bbigskip{\bigskip\bigskip}
\def\nbigskip{\bigskip\noindent}

\def\nmedskip{\medskip\noindent}
\def\buildunder#1#2{\mathrel{\mathop{\kern0pt #2}
\limits_{#1}}}
\def\buildover#1#2{\buildrel#1\over#2}

\def\qq{/\kern-.185em /}

\begin{document}

\title[Maximal Complexifications]{Maximal Complexifications
of certain  Riemannian Homogeneous Manifolds}

\author[S. Halverscheid]{S. Halverscheid}
\author[A. Iannuzzi]{A. Iannuzzi*}
\thanks{*\,Partially supported by 
University of Bologna, funds for selected research topics}

\begin{abstract}
Let $\,M=G/K \,$ be a Riemannian homogeneous manifold
with $\,dim_\C G^\C = dim_\R G\,$, where $\,G^\C\,$ denotes
the universal complexification of $\,G$.
Under certain extensibility assumptions on the geodesic
flow of $\,M$, we give a characterization of the maximal
domain of definition in
$\,TM\,$ for the adapted complex structure and show that it
is unique.
For instance, this can be done for generalized Heisenberg groups
and naturally reductive Riemannian homogeneous spaces.
As an application
it is shown that the case of generalized
Heisenberg groups yields examples of maximal domains
of definitions for the adapted complex structure which are neither
holomorphically separable, nor holomorphically 
convex.
\end{abstract}

\maketitle


\section{Introduction}
It is well known that complexifications of a real-analytic manifold
$\,M\,$ exist and 
are equivalent near $\,M,$ but differ usually very much
in nature. If a complete  real-analytic
metric on $\,M\,$ is given, one can construct canonical
complexifications of $\,M\,$ compatible with the metric
by defining an
{\it adapted} complex structure on a domain $\,\Omega\,$ of
the tangent bundle 
$\,TM$ (see \cite{GS} and \cite{LS}). 
This can be characterized by the condition that 
the ``complexification" $\,(x + iy) \mapsto y \, \gamma'(x) \in
\Omega\,$ of any geodesic
$\,x \mapsto \gamma(x)\,$ of $\,M\,$
be a complex submanifold near the zero section.
By the results of 
Guillemin-Stenzel and Lempert-Sz{\H o}ke cited above, the adapted
complex structure exists and is unique 
on a sufficiently small neighborhood of $\,M\,$.
Here $\,M\,$ is identified with its zero section in $\,TM\,$.

In particular it is natural to ask for maximal domains
around $\,M\,$ on which it exists which, by functoriality of
the definition, may be regarded as invariants of the metric, i.e.,
isometric manifolds have biholomorphic maximal domains.
For instance examples are known for
symmetric spaces of non-compact type
(\cite{BHH}), compact normal Riemannian Homogeneous
spaces (\cite{Sz2}), compact symmetric spaces (\cite{Sz1})
and spaces obtained by K\"ahlerian reduction of these
(\cite{A}). Note that in the mentioned cases
maximal domains turn out to be Stein.

The aim of this work is to characterize 
maximal domains for the adapted complex structure
for a class of Riemannian homogeneous spaces
with ``big'' isometry group.
Let $\,M=G/K\,$, with $\,G\,$ a Lie group of isometries
and $\,K\,$ compact, and
assume that $\,dim_\C G^\C = dim_\R G\,$, where $\,G^\C\,$
is the universal complexification of $\,G\,$.
Then $\,K^\C\,$ acts on $\,G^\C\,$, the left 
action on $\,M\,$ induces a natural $\,G\,$-action
on $\,TM\,$ and under certain
extensibility assumptions on the geodesic flow of
$\,M\,$ one obtains a real-analyti and  $\,G$-equivariant
map $\,P:TM \to G^\C / K^\C\,$ such that (see
Theorem \ref{sliceflow} for the precise statement)

\nmedskip
{\it The connected component of the non-singular locus
of $\,DP\,$ containing $\,M\,$ is
the unique maximal
domain on which the adapted complex structure exists.}

\nmedskip
This applies to the case of naturally reductive 
Riemannian homogeneous spaces (corollary \ref{reductive})
and of generalized Heisenberg groups (see Sect.$\,$4 and 5).

As an application it is shown that for all 
generalized Heisenberg groups such maximal 
domain is neither holomorphically separable, nor 
holomorphically convex (Proposition \ref{generalized}).
We are not aware of previous non-Stein examples.
In the case of the $3$-dimensional Heisenberg group
we determine its envelope of holomorphy
as well as a certain maximal Stein subdomain
(Proposition \ref{notstein}).
\bbigskip


\section{Preliminaries}

Here we introduce notations and briefly recall
basic results we will need in the present paper.
Let $\,M\,$ be a complete real-analytic Riemannian
manifold which will be often 
identified with the zero section in its tangent
bundle $\,TM\,$. Following 
\cite{LS} (see also \cite{GS} for an equivalent
characterization) we say that a real-analytic complex
structure defined on a domain $\,\Omega\,$ 
of $\,TM\,$ is adapted if all complex leaves of the
Riemannian foliation are submanifolds with their
natural complex structure, i. e., for any geodesic
$\,\gamma:\R \to \C \,$ the induced map $\,f: \C \to TM\,$
defined by $\,(x + iy) \mapsto y \, \gamma'(x)\,$ is holomorphic
on $\,f^{-1}(\Omega)\,$
with respect to the adapted complex structure.
Here $y \, \gamma^{\prime}(t) \in T_{\gamma(t)} M$
is the scalar multiplication in the vector space $T_{\gamma(t)} M$.

The adapted complex structure exists and
is unique on a sufficiently small neighborhood of $\,M\,$
and if $\,\Omega\,$ is a domain around $\,M\,$ in $\,TM\,$
on which it is defined,
we refer to it as an {\it adapted complexification}.
Sometimes these are called Grauert tubes.
For later use we need the following

\bigskip
\begin{lem}
\label{holomorphicity}
 Let $\,F:TM \to \C\,$ be a real-analytic map which is
         holomorphic on any complex
         leaf of $\,TM\,$ in a neighbourhood of
           $\,M\,$. Then $\,F\,$ is holomorphic
         on every adapted complexification.
\end{lem}
 
\bigskip

\begin{proof}
Following the proof of $\,$[Sz1, Prop.\,3.2 p.\,416]$\,$ 
one checks that
the restriction of $\,F\,$ to $\,M\,$  extends
to a holomorphic map $\,\hat F\,$ in a neighborhood 
$\,U\,$ of $\,M \subset TM\,$ where the adapted complex
structure $\,J_0\,$ exists and, in order to have
connected leaves, $\,U\,$ may be assumed to be starshaped.
By assumption we can also assume
that for any geodesic $\,\gamma:\R \to M\,$
the map $\,x+iy \,\mapsto \, F(y \cdot \gamma'(x))\,$ is holomorphic 
for all $\,x +iy\,$ such that $\,y \, \gamma'(x)\in U$.
Now $\,F= \hat F\,$ on 
$\,\gamma(\R)\subset M\subset TM\,$, therefore
$\,F= \hat F\,$ 
on every complex leaf, i.e., on $\,U\,$.
In particular $\,DF \circ J_0 = i \,DF\,$ on $\,U\,$
and since all maps are real-analytic the statement follows from
the identity principle.
\end{proof}

\bigskip

A real Lie group$\,G\,$ acts
on a complex manifold
$\,X\,$, i.e., $\,X\,$ is a $\,G$-manifold, if 
there exists a real-analytic surjective map
$\,G \times X \to X\,$ given by $\,(g,x) \mapsto g\cdot x\,$ such that 
 for fixed $\,g\in G\,$ the map $\,x\mapsto g\cdot x\,$
is  holomorphic  and
$\ (gh)\cdot x=g\cdot(h \cdot x) \ $ for all $\,h,g \in G\,$
and $\,x\in X$.
Furthermore if 
$\,\dim_\R \,G = \dim_\C \,G^\C\,$,
where $\,\iota:G \to G^\C\,$ is the universal
complexification of $\,G\,\,$ (see e.g. \cite{Ho}),$\,$
then $\,Lie(G^\C)= \g^\C\,$ and
one obtains an induced local holomorphic $\,G^\C$-action 
by integrating the holomorphic vector fields
given by the $\,G$-action. Here $\, \g\,$ denotes
the Lie algebra of $\,G$.

Let $\,M=G/K\,$ be a Riemannian homogeneous manifold
with $\,G\,$ a connected Lie group of isometries and 
$\,K\,$ compact, and consider the induced $\,G$-action
on $\,TM\,$ defined by $\,g \cdot w :=g_*w\,$ for all
$\,g \in G\,$ and $\,w \in TM\,$.
Then if $\,\Omega \,$ is a $\,G$-invariant adapted complexification,
as an easy consequence of the definitions 
$\,g_*\,$ is a biholomorphic extension of the isometry $\,g\,$,
i.e., $\,G \subset Aut(\Omega)\,$.

If one  assumes that
$\,\dim_\R \,G = \dim_\C \,G^\C\,$, 
then the natural map $\,\iota:G \to G^\C\,$ is an
immersion and from the universality property 
of the universal complexification $\,K^\C\,$ of $\,K\,$ it
follows that the restriction $\,\iota |_K\,$ of $\,\iota\,$ to
$\,K\,$ extends to
an immersion $\,\iota^{\C} :K^{\C} \to
G^\C\,$.
Moreover the subgroup $\,\iota^\C (K^\C)\,$ acts by
right multiplication on $\,G^\C\,$ and one has a 
commutative diagram
$$\begin{matrix}  G     \     &  \buildover{\iota}  \to \   &   \ \ G^\C \cr
                            &                            &           \cr    
              \downarrow\   &                             &     \ \downarrow \cr
                           &                               &                   \cr
         G/K     \   &       \to \    & \  \ \  \ G^\C/\iota^\C(K^\C) \,.\cr 
\end{matrix} $$

\nsmallskip
Also notice the $\,G$-action on $\,G^\C/\iota^\C(K^\C)\,$ 
defined by $\,g \cdot h\,\iota^\C(K^\C) := \iota(g)\, h\,\iota^\C(K^\C)\,$
for all $\,g \in G\,$ and $\,h \in G^\C\,$. 
\medskip
\begin{lem}
\label{goodquotient}
Let $\,G\,$ be a connected Lie group, $\,K\,$ a compact 
subgroup and assume that
$\dim_\C G^\C = \dim_\R G$.
Then 
$G^\C/\iota^\C(K^\C)\,$ is a complex $\,G$-manifold and 
$\,\dim_\C \,G^\C/\iota^\C(K^\C) = \dim_\R\, G/K.$
\end{lem} 

\bigskip
\begin{proof}
One needs to show that $\,\iota^\C(K^\C)\,$ is closed in
$\,G^\C\,$.
Note that $\,G^\C/\iota^\C(K^\C)\,$ is the orbit
space with respect to the $\,K^{\mathbb{C}}$-action on $G^{\mathbb{C}}\,$
defined by $\,k \cdot h := h\,\iota(k^{-1})\,$ for all
$\,k \in K^\C\,$ and $\,h \in G^\C\,$.
Since $\,G^\C\,$ is Stein (\cite{He}) and 
$\,K^\C\,$ is reductive it follows that
every fiber of the categorical
quotient  $\,G^{\mathbb{C}} \to G^{\mathbb{C}} \qq K^{\mathbb{C}}\,$
is equivariantly biholomorphic to an affine algebraic variety on which 
$\,K^{\mathbb{C}}\,$ acts algebraically (\cite{Sn}). In particular
there exists at least one closed $\,K^\C$-orbit and consequently 
$\,\iota^\C (K^{\mathbb{C}})\,$ is closed in  $\,G^{\mathbb{C}}\,$. 
Thus  $\,G^{\mathbb{C}}/\iota^\C(K^{\mathbb{C}})\,$ is a complex 
$\,G$-manifold and by construction its complex dimension is
$\,\dim_\R\,G/K \,$. 
\end{proof}

\bigskip


\section{A characterization of maximal adapted complexifications}
\bigskip

If $\,M=G/K\,$ is a symmetric space of the non-compact
type, then $\,G^\C/K^\C\,$ is a natural candidate for a
complexification of $\,M\,$ and
there exists a $\,G$-equivariant map
$\,P:TM \to \,G^\C/K^\C\,$
embedding 
holomorphically a maximal
adapted complexification of $\,M\,$ (see \cite{BHH}, \cite{Ha}, \cite{AG}).
As a matter of
fact one may show that $\,DP\,$ is singular on 
the border $\,\partial \Omega\,$ of $\,\Omega$.

Here we consider a Riemannian homogeneous manifold
$\,M=G/K\,$
endowed with the additional data of a certain
real-analytic 
$\,G$-equivariant map
$\,P\,$ from $\,TM\,$ to a suitable 
complex $\,G$-manifold,
characterizing a maximal adapted complexification
$\,\Omega_M\,$ as the connected component of
$\,\{\,DP$ not singular $\}\,$ containing $\,M$.
Unicity of $\,\Omega_M\,$ follows.

The existence of such a data is proved when
$\dim_\C G^\C = \dim_\R G$ and the geodesic flow ``extends''
holomorphically on $G^\C/\iota^\C(K^\C)\,$ 
(cf. Lemma\,\ref{goodquotient}).
As a consequence the characterization applies to the case 
of naturally reductive Riemannian homogeneous spaces 
and of generalized Heisenberg groups.

\bigskip
\begin{prop}
\label{equivariant-theorem}
Let $\,M=G/K\,$ be an $n$-dimensional Riemannian homogeneous
space and $\,X\,$ a $\,G$-complex
manifold of complex dimension $n$
such that the induced local $\,G^\C\,$-action is
locally transitive. Assume there exists a real-analytic
map $\,P:TM \to X\,$ which is 
\begin{align*}
 \ \ &  i)  \ \ \ {G-equivariant \ and}   &   \\
 \ \ &  ii)  \ \ 
{holomorphic \ on\  every \ complex \ leaf \ of \ TM.}  & 
\quad \quad \quad \quad \quad \quad \quad
\quad \quad \quad  \quad 
\end{align*} 

\medskip
\noindent
Then the connected component $\,\Omega_M\,$ of 
$\, \{ \,p \in TM \, : \, DP_p {\rm \ \ not\  singular\ }\,\}\,$
containing $\,M\,$ is the unique maximal adapted complexification
and $\,P|_{\Omega_M}\,$ is locally
biholomorphic. 
\end{prop}

\bigskip

\begin{proof} 
First we show that $\,\Omega_M\,$ is well defined, i.e.,
$\,DP\,$ has maximal rank along $\,M\,$. 
Since from Lemma \ref{holomorphicity} it follows that
$\,P\,$ is holomorphic on $\,M\,$ with respect to the
adapted complex structure, this is
a consequence of the following

\nmedskip
{\it Claim:} Assume that $\,P\,$ is holomorphic in
$\,p \in TM$. Then $\,DP_p\,$ has maximal rank.

\nmedskip
{\it Proof of the claim:} $\ \ $
Since $\,G^\C\,$ acts locally transitively
on $\,X\,$, there exist elements $\,\xi_1, \cdots, \xi_n\,$
of $\, \mathfrak g \,$ such that the induced 
vector fields $\,\xi_{X,1}, \cdots,
\xi_{X,n}\,$ on $\,X\,$
span a totally real and maximal dimensional
subspace $\,V_{P(p)} \,$ 
of  $\,T_{P(p)}X\,$,
where
$$\,\xi_{X,j}(x):=\frac{d}{dt}\biggr|_0 \exp_{G^\C}(t\xi_j) \cdot x\,$$
for $\,j=1, \cdots, n\,$ and all $x \in X$. By equivariance it follows
that $\,DP_p(V_p) =V_{P(p)}\,$, where  $\,V_p\,$ is the 
subspace of
$\,T_pTM\,$ spanned by 
$\,\xi_{TM,1}, \cdots , \xi_{TM,n}\,$, here 
$\,\xi_{TM,j}(q):=\frac{d}{dt}\biggr|_0 \exp_{G}(t\xi) \cdot q\,$
for all $\,q \in TM$.
In particular $\,{\mathrm dim}_\R V_p = n\,$ and since 
$\,P\,$ is holomorphic in $\,p\,$, 
$\,V_p\,$ is totally real and  
$\,DP_p\,$ has maximal rank, proving the claim.

\medskip
Now we see that the pulled-back complex structure $\,J_o\,$
on $\,\Omega_M\,$ of the complex structure $\,J\,$ on $\,X\,$
is the adapted complex structure.
For this consider a complex leaf $\,f:\C \to TM\,$ defined by
$\,f(x+iy):=y \cdot \gamma'(x)\,$, where $\,\gamma\,$ is a 
geodesic of $\,M\,$, and note that by $ii)$
$$\,DP \circ Df(i \eta) =DP \circ J_o \circ
Df(\eta)\,$$ for all $\,\eta \,$ tangent in 
$\,f^{-1}(\Omega_M)\,$. Since $\,DP\,$
has maximal rank on $\,\Omega_M\,$, then
$$\,Df(i \eta) = J_o \circ Df(\eta)\,$$
showing that $\,J_o\,$ is the
adapted complex structure. 
In particular $\,P|_{\Omega_M}\,$ is locally
biholomorphic.

In order to prove maximality,
assume that $\,J_o\,$ extends analytically in
a neighborhood of a certain $\,p\in \partial \Omega_M
\subset TM\,$. By construction $DP \circ J_o = J \circ DP \,$
on $\,\Omega_M\,$
and since all maps are real-analytic
$\,P\,$ is holomorphic in $\,p\,$.
Then the above claim shows that 
$\,DP_p\,$ has maximal rank, contradicting the definition
of $\Omega_M$.

Finally we want to show that any adapted complexification $\,\Omega\,$
is contained in $\,\Omega_M$. If this is not the case,
there exists a point $\,p\,$ in $\,\Omega \cap \partial \Omega_M$
and from Lemma \ref{holomorphicity} it follows that
$\,P|_\Omega\,$ is holomorphic. In particular $\,P\,$
is holomorphic in $\,p\,$ and one obtains a contradiction
arguing as above. Thus $\,\Omega_M$ is unique and this concludes
the statement.
\end{proof}

\bigskip

Now we determine a class of Riemannian homogeneous spaces
to which Proposition \ref{equivariant-theorem} 
may be applied in order to determine
the maximal adapted complexification.

\bigskip
\begin{theorem}
\label{sliceflow}
$\ \ \ \ $ Let $\,M =G/K\,$ be a Riemannian homogeneous
space with
$\, \dim_\R \,G = \dim_\C \,G^\C\,$
and assume  there exists a map
$\, \varphi: \C \times T_K M \mapsto \mathfrak g^\C \,$ 
real-analytic and holomophic on the first component
such that $\,\varphi(\R \times T_K M) \subset \mathfrak g\,$
and $\, t \mapsto \exp_G \circ \varphi (t,v)K\,$ is the unique geodesic
tangent to $\,v\,$ at $\,0\,$ for all $\,v \in T_K M\,$.
Then the map 
$$P: TM  \mapsto G^\C / \iota^\C(K^\C)$$
defined by
$$P(g_*(v)):= \iota(g)\,\exp_{G^\C}(\varphi(i,v))\,
\iota^\C(K^\C)$$ 

\medskip
\noindent
for all $\, g \in G\,$ and $\,v \in T_KM$ is as in
Proposition 
\ref{equivariant-theorem}. In particular
the connected component $\,\Omega_M\,$ of 
$\,\{ \,p \in TM \, : \, DP_p {\rm \  \ not\  singular\ }\,\}\,$
containing $\,M\,$
is the maximal adapted complexification and $\,P|_{\Omega_M}\,$ is locally
biholomorphic.
\end{theorem}

\bigskip
\begin{proof} 
In order to prove that $\,P\,$ is well defined we need to show
that if $\,w=k_*v\,$ for some $\,k \in K\,$ and $\,v \in T_K M\,$
then $\,P(w)=P(k_*v)\,$, i. e., 
\begin{equation}
\label{welldefined}
\, \exp_{G^\C}(\varphi(i,w))\,\iota^\C(K^\C) =
\iota(k)\,\exp_{G^\C}(\varphi(i,v))\,\iota^\C(K^\C)\, .
\end{equation}

\nsmallskip
For this note that $\, t \mapsto k\,\exp_G \circ \varphi (t,v)K\,$
is the unique geodesic tangent to $\, w \,$ at $\,0\,$ in $K$,
thus 
$$\,\exp_G \circ \varphi (t,w)K\,=k\,\exp_G \circ \varphi (t,v)K\,.$$

\nmedskip
Then the commutativity of the diagram 
$$\begin{matrix}  \g \ \ \ \ \ \   &   \to \  &   \g \ ^\C \ \ \ \  \cr
                            &                            &           \cr
         \downarrow\, \exp_G   &          &     \ \downarrow \, \exp_{G^\C}  \cr
                           &                               &                   \cr
         G \ \ \ \ \ \ \  & \buildover{\iota}   \to  \    & G^\C \ \ \ \ \ \cr
\end{matrix} $$

\noindent 
implies that 
$$\,\exp_{G^\C} \circ \varphi (t,w)\,\iota^\C(K^\C)\,=
\iota(k)(\,\exp_{G^\C} \circ \varphi (t,v)\,) \iota^\C(K^\C)\,$$

\nmedskip
for all $\,t \in \R\,$
and equation (\ref{welldefined}) is a consequence of the identity principle
for holomorphic maps.

 Now define $\,\varPhi(z,v):=\exp_{G^\C} \circ \varphi (z,v)\,$
 for all $\,z \in \C\,$ and $\,v \in T_K M\,$
and, in order to simplify notations, assume that 
the canonical immersion $\,\iota :G \to G^\C\,$ is injective 
 so that once we identify $\, G \,$ with  $\iota(G)$ the curve
 $\, \gamma(t) := \varPhi(t,v)K\,$ is the unique  geodesic tangent to 
$\, v\,$ at $\,0\,$.
In what follows it is easy to check that all
arguments apply to the case where $\,\iota\,$ is
a non-injective immersion.

Fix $\,x \in \R\,$, let $g:=\varPhi(x,v)\,$ and note that 
$$ y \mapsto g\,\varPhi(y, g_*^{-1}\gamma'(x))K$$

\nmedskip
is the unique geodesic tangent to $\,\gamma'(x)\,$
at $\,0\,$. Therefore one has 
$$\varPhi(x,v)\,\varPhi(y,g_*^{-1}\gamma'(x))K =
\varPhi(x+y,v)K$$

\nmedskip
for all $\,y\in \R\,$ and  by the identity principle it follows that 
\begin{equation}
\label{flowrule}
 \varPhi(x,v)\,\varPhi(z,g_*^{-1}\gamma'(x))\,\iota^\C(K^\C) =
\varPhi(x+z,v)\,\iota^\C(K^\C)
\end{equation}

\nsmallskip
for all $\,v \in T_K M, \ x \in \R\,$ and $\, z \in \C\,$.
For $\, h \in G \,$ and $\,v\in T_KM\,$ consider the unique 
geodesic  $\,\tilde \gamma(t) :=h \cdot \gamma(t)\,$
tangent to 
$\,h_*(v)\,$ at $\,0\,$. One has 
$$P(y\tilde \gamma'(x))=P(h_*y \gamma'(x))=
h\,P(g_*g_*^{-1}y \gamma'(x))=
hg\,\varPhi(i,y \,g_*^{-1}\gamma'(x))\,\iota^\C(K^\C)=$$
$$h\,\varPhi(x,v)\,\varPhi(iy,\,g_*^{-1}\gamma'(x))\,\iota^\C(K^\C)=
h\,\varPhi(x +iy,v)\,\iota^\C(K^\C),$$

\medskip
\noindent
where we used (\ref{flowrule}) and the fact that 
$\, \varPhi(z,yv)\,\iota^\C(K^\C)= \varPhi(zy,v)\,\iota^\C(K^\C)\,$
for all $\,z \in \C\,$, since this holds for all
$\,z \in \R\,$.
As a consequence the map $\, (x+iy) \mapsto P(y\tilde \gamma'(x))\,$
is holomorphic for all geodesics 
$\, \tilde \gamma \,$ of $\,M\,$, i.e., $\,P\,$
is holomorphic on every complex leaf of 
$\,TM$.

Finally the map $\, P \,$ is $\,G$-equivariant by construction
and the $\,G$-action on $G^\C/\iota^\C(K^\C)$ induces
a holomorphic $G^\C$-action which may be 
obtained through left multiplication on $G^\C$.
Thus it is obviously transitive and this yields the 
statement.
\end{proof}                  

\bigskip


Now let $\,M\,$ be a naturally 
reductive Riemannian homogeneous space
and $\,M=G/K\,$ be a natural realization of $\,M\,$, i.e., there exists
a reductive decomposition $\,\g = Lie(K) \oplus 
\mathfrak{m}\,$ of the Lie algebra of $\,G\,$ such that
every geodesic in $\,M\,$ is the orbit of
a one parameter subgroup of $\,G\,$ generated by an element of
$\,\mathfrak{m}\,$ (see e.g. \cite{BTV}). 
Consider the natural projection 
$\, \Pi: G \to M\, $ and note that $\,D\Pi_e(\mathfrak{m})  =
T_K\,M$, where $\,e\,$ is the neutral element of $\,G$. Denote by $\,L: T_K M \to \mathfrak{m}\,$ the
inverse of the restriction of 
$\,D\Pi_e\,$ to $\, \mathfrak{m}\,$. Since 
$\,L \,$ is linear it extends $\,\C$-linearly from 
$\,(T_K M)^\C\,$ to  $\, \mathfrak{m}^\C\,$
and the map $\, \varphi: \C \times T_K M \to 
\mathfrak{g}^\C\,$ defined by $\, \varphi(z,v):=zL(v)\,$
is as in the above Theorem. Therefore one has
\bigskip


\begin{cor}
\label{reductive}
Let $M =G/K$ be a natural realization of a naturally
reductive Riemannian homogeneous space and assume that
$\dim_{\mathbb{R}}G = \dim_{\mathbb{C}}G^{\mathbb{C}}$.
Then the map 
$$P: TM  \to G^{\mathbb{C}}/\iota(K^{\mathbb{C}})$$
defined by
$$P(g_*(v)):= \iota(g)\,\exp_{G^\C}(iL(v))\,\iota(K^\C)$$ 

\nmedskip
for all $\, g \in G\,$ and $\,v \in T_KM$ meets the conditions of
Proposition \ref{equivariant-theorem}.
In particular
the connected component $\,\Omega_M\,$ of 
$\,\{ \,p \in TM \, : \, DP_p {\rm \ \ not\  singular\ }\,\}\,$
containing $\,M\,$ is the maximal adapted complexification
and $\,P|_{\Omega_M}\,$ is locally
biholomorphic.
\end{cor}

\bbigskip


\section{The 3-dimensional Heisenberg group}

\bigskip

Here we apply results of the previous section in order to
give a concrete description of the unique maximal adapted complexification
for the 3-dimensional Heisenberg group. It
turns out that such domain is neither holomorphically
separable, nor holomorphically convex.
We also determine its envelope of holomorphy and
a particular maximal Stein subdomain.
We remark that in all previous examples we are aware of,
maximal adapted complexifications are Stein.

\nbigskip
Consider the 3-dimensional Heisenberg group defined as
a subgroup of $\, GL_3(\R)\, $ by

{\tiny
\begin{align*}  H := \left \{ \left( \begin{array}{ccc}
                                   1 & \alpha   & \gamma \\
                                   0 &   1      & \beta   \\
                                   0 &   0      &    1
                \end{array} \right)\ : \ 
\alpha, \beta, \gamma \in \R \right \}, 
& 
\end{align*} 
}

\nmedskip
fix the inner product of the tangent space $\,T_eH\,$
in the neutral element $\,e\,$ 
for which the canonical basis determined by 
the global natural chart $\,(\alpha, \beta, \gamma)\,$
is orthonormal and let 
$\,(a,b,c)\,$ be coordinates of $\,T_eH\,$ with respect
to this basis. 
Endow  $\,H\,$ with the induced
$\,H$-invariant metric
$$ (\mathrm{d} \alpha)^2 + (\mathrm{d} \beta)^2 + 
   (\mathrm{d} \gamma - \alpha \cdot \mathrm{d} \beta)^2,$$ 

\nsmallskip
let $\,\h=Lie(H)\,$ and define
$\,\varphi:\R \times T_eH \,\to\, \h \,$ by

{\tiny
\begin{align*}
\varphi(t,(a,b,c)) \,:=\,
 \left (\,a \frac{\sin(tc)}{c} - b \frac{1 -\cos(tc)}{c}\ , \ 
       b \frac{\sin(tc)}{c} + a \frac{1 -\cos(tc)}{c}\ , \ 
        \left ( t +\frac{a^2 +b^2}{2c^2} \left ( t - \frac{\sin(tc)}{c}
        \right ) \right) c \,
\right ), & 
\end{align*} 
}

\noindent
where the coordinates of $\,\h\,$ are induced by those of $\,T_eH\,$
via the natural identification $\,\h \cong T_eH$.
Note that all singularities are removable and
consequently 
$\,\varphi \,$ is real-analytic. Following [BTV, Th.\,p.\,31]
one checks that
$\,t \mapsto \exp_H \circ\varphi(t,(a,b,c))\,$ is the unique
geodesic tangent to $\,(a,b,c)\,$ at $ \,0\,$. Furthermore by expanding
the power series it is easy to verify  that
$\, \varphi(\, \cdot \,  ,(a,b,c))\, $
extends holomorphically on $\,\C\,$ to 
$\,(T_eH)^\C\,$ and by considering the polar decomposition
$\,H \times \h \to H^\C\,$ of $\,H^\C\,$ 
given by
$\,(g,\xi) \,\mapsto \,g\,\exp_{H^\C}(i\xi)\,$,
one obtains real-analytic functions 
$\, (a,b,c) \mapsto h_{(a,b,c)} \in H\,$ and 
$\,(a,b,c) \to \xi_{(a,b,c)} \in \h\,$ such that
$$
\exp_{H^\C} \circ \varphi(i,(a,b,c))\,=\,
          h_{(a,b,c)}\, \exp_{H^\C} (i\xi_{(a,b,c)}).
$$

\nmedskip
Define $\,P:TH \,\to\, H^\C \, \cong \,H \times \h \,$ by
$$\,g_*(a,b,c) \, \mapsto \, g\, \exp_{H^\C}\circ \varphi(i,(a,b,c)) \,\cong
\,\left ( \, gh_{(a,b,c)}\, , \, \xi_{(a,b,c)} \,\right ) .$$


\nmedskip
Then Theorem \ref{sliceflow} implies that the connected
component $\,\Omega_H\,$
containing $\,H\,$  of $\,\{\,DP$ not singular $\}$ 
is the maximal adapted complexification.
Note that since $\,P\,$ is  $\,H$-equivariant, 
$\,\Omega_H\,$ is $\,H$-invariant.
Moreover $\,T_eH\, $
is a global slice for the $\,H$-action on $\,TH\,$, i.e., the map
$\,H \times T_eH \to TH\,$ given by $\,(g,(a,b,c))\to g_*(a,b,c)\,$
is a $\,H$-equivariant real-analytic diffeomorphism, thus 
 $\,\Omega_H\,$ is completely determined by its slice $\,
\Omega_H\cap T_eH$.

Furthermore $\,H\,$ acts 
freely on the first component of  $\,H \times \h \,$,
then $\,H$-equivariance of $\, P \,$  implies
that $\, DP_{g_*(a,b,c)} \,$ has maximal rank if and only if
$\, D\tilde P_{(a,b,c)} \,$ has maximal rank,
where
$\, \tilde P\, := \,p_2 \circ P|_{T_eH}:T_eH \,\to \,\h
\,$ is given
by 

{\tiny
\begin{align*}
\tilde P(a,b,c) \,=\, \xi_{(a,b,c)}\, =\,
\left ( \,a \frac{\sinh(c)}{c}\ , \ 
       b \frac{\sinh(c)}{c} \  , \ 
        \left ( 1 +\frac{a^2 +b^2}{2c^3} \left ( c- \sinh(c)\,\cosh(c)
\right ) \right) c \,\right ).
& 
\end{align*} }

\nmedskip
Here $\,p_2:H \times \h \,\to \,\h \,$ is the canonical
projection.
It follows that $\,\Omega_H=H\cdot O_0\,$, where
$\,O_0\,$ is 
 the connected component of 
$\{ \,\det(D\tilde P)\, \not= \,0 \}$ containing  $\,0\,$ in
$\,T_eH\, $.
Now a straightforward computation shows
that 

{\tiny
\begin{equation*}
\det(D\tilde P_{(a,b,c)} )\,=\, \frac{\sinh(c)}{c} \left (\,\frac{\sinh(c)}{c}
+ (a^2 +b^2) \left (\frac{\sinh(c) -c \,\cosh(c)}{c^3} \right ) \right),
\end{equation*} }

\noindent
therefore

{\tiny
\begin{align*}
O_0= \left \{ \,(a,b,c) \in T_eH \ : \  a^2 + b^2 < \frac{ c^2 \, 
      \sinh(c)}{c\,\cosh(c) - \sinh(c)} \right \}.
& 
\end{align*} }

\nmedskip
We want to discuss injectivity of 
$\,P|_{\Omega_H}:\Omega_H \, \to \, \, H^\C \,
\cong \,H \times \h \,$ and again this is equivalent
to injectivity of $\,\tilde P|_{O_0}\,$.

Note that
$\,\tilde P\,$ is equivariant with respect to rotations
around the $\,c$-axis as well as to the 
reflection $\, \sigma \,$ with respect to the  plane
$\,\{ \, c=0 \, \}\,$.
In particular for any 
{\tiny
\begin{align*}
\,  (a,b,c) \in \left \{ 1 +\frac{a^2 +b^2}{2c^3}
\left ( c - \sinh(c)\,\cosh(c)\right )  \,=\,0  \right \}\,
& 
\end{align*} }

\noindent
one has $\, \tilde P(a,b,c) =\tilde P(a,b,-c) =
\left (\,a\,\frac{\sinh(c)}{c}, b\,\frac{\sinh(c)}{c}, 0 \,\right ) \,$. 
Therefore we are induced to investigate the
domain 

{\tiny
\begin{align*}
O_1:= \left \{ \,(a,b,c) \in T_eH \ : \  a^2 + b^2 
< \frac{ 2c^3}{\sinh(c)\cosh(c) - c} \, \right \}.
& 
\end{align*}}

\medskip
\begin{lem}
\label{injectivdomain}
The domain $O_1\,$ is the  maximal $\,\sigma \,$-invariant
subdomain of $\,O_0\,$
containing $\,0\,$ on which 
$\, \tilde P\,$ is injective. In particular
$\, \tilde P|_{O_0}\,$ is not injective.
\end{lem} 

\bigskip
\begin{proof}
Let $\, f_j\,$ be the real function defining
$$\,O_j\,=\,\{ \,(a,b,c) \in T_eH \ : \  a^2 + b^2 
< f_j(c)\,\}\,$$ for $\,j=0,1 \,$. 
First we want to show that $\,O_1\,$ is a
subdomain of $\,O_0\,$, i.e, 
$$ f_1(c)\,\leq\,f_0(c)$$

\nmedskip
for all $\,c \in \R\,$, which is equivalent to 
$$2\cosh(c) \leq \frac{\sinh(c)}{c} +\frac{\sinh^2(c)}{c^2}\cosh(c).$$

\nmedskip
Expanding in power series one obtains
{\tiny \begin{align*}
2 \,+\, c^2 &+\frac{1}{12}c^4\, + \, \cdots  \ \,\leq\,\, 
 \left( 1\,+\frac{1}{6}c^2\,+\,\frac{1}{120}c^4\,+\,\cdots \right )
\,+\,\left ( 1\, +\,\frac{5}{6}c^2\,+\,\left(\frac{5}{6}\,+
\frac{1}{24}\,+\,\cdots\, \right)c^4 \,+ \, \cdots \,\right ). & 
\end{align*}}

\noindent
All coefficients
are non-negative and one easily checks that for $\, k \geq 2\,$
the
coefficient of $ \,c^{2k}\,$ in the last series on
the right side is strictly greater than that in the series on the left,
hence $\,O_1\subset O_0\,$. Moreover 
$\,\partial O_1 \cap \partial O_0\,= \, \{\,
(a,b,0) \in T_eH \ : \ a^2 + b^2 =3 \, \}$, thus
$\,O_1\,$ is a proper subdomain of $\,O_0$.

Furthermore by the previous remarks any 
$\,\sigma$-invariant domain containing 
$\,0\,$ on which $\, \tilde P\,$ is injective
is necessarily contained in $\,O_1 \,$.

Assume that there exist $\,(a',b',c'), \,(a'',b'',c'') \in O_1\,$
such that $\, \tilde P(a',b',c')=\tilde P(a'',b'',c'')=:(A,B,C).$
If $\,C=0\,$ then $\,c' =c''=0\,$ and consequently
$\,a'=a'' =A   \,$ and $\,b=b''=B\,$.
If $\,C\not=0\,$ by eventually acting with
$\,\sigma\,$ and a rotation 
around the $\,c$-axis we may assume
that $\, a,A\geq 0\,$, $\,b=B=0\,$ and $\, c >0\,$ .
Now one has 
$$\,a'\frac{\sinh(c')}{c'}=a''\frac{\sinh(c'')}{c''}=A\,$$
therefore
$\,(a',0,c')\,$ and $ \,(a'',0,c'')\,$ lie on the same 
level curve
$\,\rho_A:\R \to T_eH\,$ given by 
$$\rho_A(t):=\left(\,A \frac{t}{\sinh(t)},0,t\, \right ).$$
One has the following 

\nmedskip
{\it Claim:} $\ \ $Let $\, A\geq 0\,$ and $\, t_0\in\R^{\geq 0} \,$
such that $\rho_A(t_0) \in \overline O_0\,$. Then 
$\rho_A(t) \in O_0\,$ for all $\, t >t_0\,$.

\nmedskip
{\it Proof of the claim:} $\ \ $
One needs to show that 
$A^2\frac{t^2}{\sinh^2(t)} \,< \, f_0(t)$ for all $\,t>t_0\,$,
that is
{\begin{equation}
\label{disugua}
A^2\,<\,\frac{\sinh^2(t)}{t^2} f_0(t).
\end{equation}}

\noindent
By expanding
in power series as above one
has the estimate
$$
2t\,\cosh^2(t) -3\cosh(t)\sinh(t) + t >0,
$$

\nmedskip
for all $\,t>0\,$,
which by a straightforward computation
implies that the derivative of the
function at the right hand side of (\ref{disugua}) is positive for
all $\,t >0\,$, proving the claim.

\nmedskip
Now let $\,t_0 := \min(c',c'')\,$ and note that 
since $\,O_1 \subset O_0\,$, then as a consequence of
the above claim
there exists $\,\epsilon >0\,$ such that
$\,\rho_A(t) \in O_0 \,$ for $\,t >t_0-\epsilon$. 
In particular $\,(a',0,c')\,$ and $ \,(a'',0,c'')\,$ lie
in the same connected real one dimensional 
submanifold $\,N:=\rho_A(t_0-\epsilon, \infty)\,$
of $\,O_0\,$
and $\,\tilde P|_{N}:N \,\to\, \{(A,0,\,\cdot\,) \in T_eH\} \, \cong \,\R\,$
is locally diffeomorphic. Then a classical
argument implies   that $\,\tilde P|_{N}\,$
is injective, thus $\,(a',0,c')\,= \,(a'',0,c'')\,$
as wished.
\end{proof}

\nbigskip
We also  want to determine the image  
of $\,P|_{\Omega_H}\,$ in   
$\,H^\C\,$. Note that $\,P(\Omega_H)\,$ is $\,H$-invariant and
the polar decomposition implies  that 
$\,\exp_{H^\C}(i\g)\,$ is a global slice for the $\,H$-action
on $\,H^\C\,$. Then this can be achieved by describing
 $\,exp_{H^\C}(i\tilde P(O_0))=P(\Omega_H)\cap \exp_{H^\C}(i\g)\,$.

\bigskip
\begin{lem}
\label{imagedomain}$\,\tilde P(O_0)
\,=\, \h \setminus \{ \,(A,B,C) \in \h
\ :\  A^2 +B^2 =3, \ C=0 \}.$
\end{lem} 
\bigskip
 
\begin{proof}
Let  $\,(a,0,c) \in \{a^2 =f_1(c)\} \subset \partial O_1\,$  with 
$ \, c  >  0\,$.
From the proof of Lemma 
\ref{injectivdomain} it follows that 
$\,(a,0,c) \in O_0\,$. Since 
$\,\tilde P(a,0,c)\,=\,(\,a\frac{\sinh(c)}{c},0,0\,)\,$
and
$$\sqrt{f_1(0)}  =  \sqrt 3\quad {\rm and} \quad \quad \ \ 
\,\lim_{c \to\infty} \sqrt{f_1(c)} \, \frac{\sinh(c)}{c} = \infty\ $$ 
\nmedskip
it follows that  $\,(A,0,0) \in \tilde P(O_0)\,$ for all 
$\,A >\sqrt 3\,$.

For $\,A> \sqrt 3 \,$  let  $\,(a,0,c) \in O_0\,$ 
such that 
$\,\tilde P(a,0,c)\,=\,(\,A,0,0\,)\,$.
By the claim
in Lemma \ref{injectivdomain}
one has
$\, \rho_A(t) \in O_0\,$ for all $\, t \geq c\,$.
Moreover one sees that 

\begin{equation}
\label{limit}
\,\tilde P(\rho_A(t))=(A,0,C_A(t)) \quad \quad {\rm with} \quad \quad
\,\lim_{t \to\infty}C_A(t) = \infty\,.
\end{equation} 

\nmedskip
Then by $\,\sigma$-invariance
of $\,O_0 \,$ and $\,\sigma$-equivariance
of $\,\tilde P \,$ it follows that $\,(A,0,C) \in
\tilde P(O_0)\,$ for all $\,C\in \R,$
and $\,A >\sqrt 3\,$.

Now note that $\,\tilde P(a,0,0)=(a,0,0)\,$ and 
$\,f_0(0)=3\,$, thus $\,(A,0,0) \in \tilde P(O_0)\,$
for all $\,A< \sqrt 3\,$ and arguing as above it
follows that $\,(A,0,C) \in \tilde P(O_0)\,$
for all $\,C\in \R\,$ and $\,A < \sqrt 3\,$.

Finally $\, \rho_{\sqrt 3}(0) \in \partial O_0\,$,
thus  $\, \rho_{\sqrt 3}(t) \in O_0\,$ for all
$\,t>0\,$. It follows that
$\,t \mapsto \tilde P(\rho_{\sqrt 3}(t)) \,$
is injective for $\,t>0\,$
and since 
$$\,\lim_{t \to 0^+}C_{\sqrt 3}(t) = 0\
\quad {\rm and} \quad
\,\lim_{t \to\infty}C_{\sqrt 3}(t) = \infty\,$$

\noindent
then   
$\,(\sqrt 3,0,C) \in \tilde P(O_0)\,$
if and only if $\,C\not=0\,$ .

The statement follows from the
invariance of $\,O_0 \,$ and the equivariance
of $\,\tilde P \,$ with respect to the group
of rotations around the  $\,c$-axis.
\end{proof}


In the picture below one sees the border
of $\,O_0\,$ and $\,O_1\,$ determined by
$\,f_0 \,$ and $\, f_1\,$ respectively
as well as the level curve $\,\rho_2\,$
in the upper half-plane $\,\{\, b =0, \ a \geq 0 \}$
of $\,T_eH$.
Since $\,O_0\,$ and $\,O_1\,$ are invariant
with respect to rotations around the $\,c$-axis
this completely determine their shape and,
by $\,H$-equivariance, that of $\,\Omega_H$.

\begin{figure}
  \begin{center}
    \leavevmode 
    \epsfxsize=0.5\textwidth     
    \epsffile{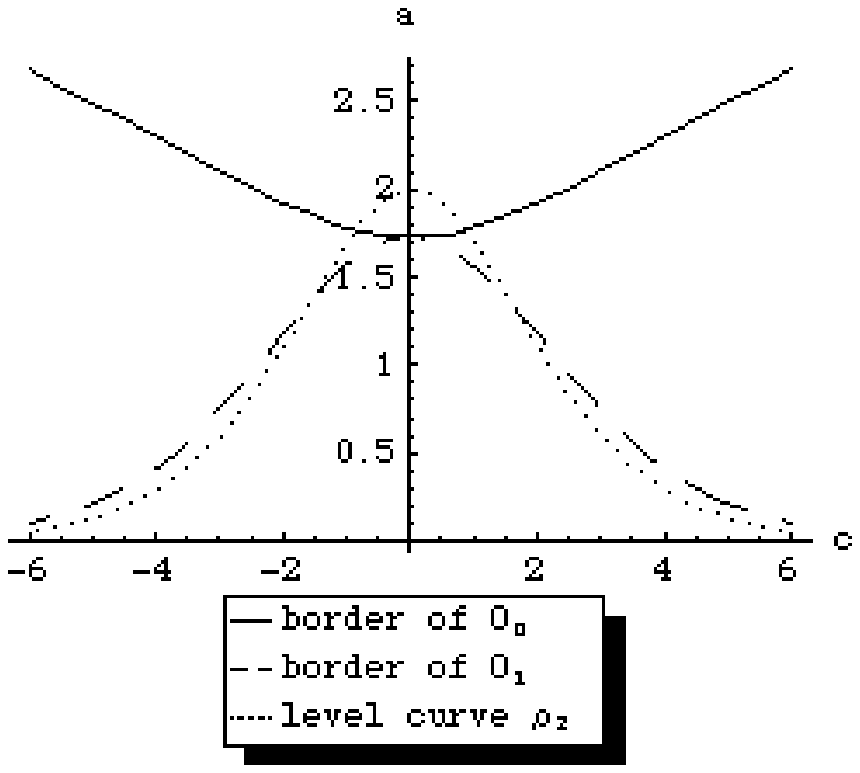}
  \end{center}
\end{figure} 


\bigskip
\begin{prop}
\label{notstein}
The maximal domain
$\Omega_H\,$ is neither holomorphically separable,
nor holomorphically
convex, its envelope of holomorphy
is biholomorphic to
$\,\C^3$. The maximal $\,\sigma$-invariant
Stein subdomain of $\,\Omega_H\,$ is
biholomorphic to $\,\{\, (z_1,z_2,z_3) \in 
\C^3 \ : \ (\mathrm{Im}\,z_1)^2\, +\, (\mathrm{Im} z_2)^2
 < 3 \, \}$.
\end{prop}

\medskip
\begin{proof}
From [CIT, Prop.\,p.\,543] it follows that
any holomorphically separable 
Riemann $\,H$-domain over $\,H^\C\,$ is 
univalent. Moreover Lemma \ref{injectivdomain} implies that
$\,P|_{\Omega_H}:\Omega_H \, \to \, \, H^\C \,$
is not injective, therefore $\Omega_H\,$ is not
holomorphically separable.

By a result of Loeb ([L,\,Th.\,p.\,186]),
a Stein $\,H$-invariant domain
$\,U\,$
of $\,H^\C\,$ is ``geodesically" convex, i.e.,
it is convex with respect to all curves of the form
$\, t \mapsto g\,\exp_{H^\C}(it \xi)\,$, with 
$\,g \in U\,$ and $\,\xi \in \h\,$.
Since $\,H^\C\,$ admits polar decomposition
and $\,U\,$ is $\,H$-invariant
it is enough to consider curves of the form
$\,\exp_{H^\C}(i \eta)
\,\exp_{H^\C}(it \xi)\,$, with 
$\,\exp_{H^\C}(i \eta) \in U\,$ and
$\,\xi \in \h\,$.
Furthermore for a two steps nilpotent Lie group one has
$$\,\exp_{H^\C}(i\eta)\,\exp_{H^\C}(it\xi)\,
=\,\exp_{H^\C}( i\eta \,+ \,it\xi \,- \,\frac{t}{2}[\eta, \xi] )\,=$$
$$\,\exp_{H^\C}( - \,\frac{t}{2}[\eta, \xi] )\,
\exp_{H^\C}( i(\eta \,+ \, t\xi)\,)$$

\nmedskip
and using $\,H$-invariance one more time
we conclude that if $\,U=H \cdot \,
\exp_{H^\C}(iD)\,$, with  $\,D \,$ a domain in 
$\,\h\,$, is Stein then  $\,D \,$ is 
convex in the usual affine sense.

Since $\,P(\Omega_H)= H\cdot exp_{H^\C}(i\tilde P(O_0))\,$
and as a consequence of Lemma \ref{imagedomain} the domain
$\,\tilde P(O_0)\,$ is not convex,
then $\,P(\Omega_H) \,$ is not Stein.
Now $\,H^\C\,$ is Stein and $\,P|_{\Omega_H}\,$ is locally
biholomorphic, therefore  
by \cite{R} there exists a commutative diagram
$$
\begin{matrix}
\Omega_H   & \buildover j \longrightarrow   &  \hat \Omega_H\cr
                    &                                                        &               \cr 
                    &  P \searrow                              &\downarrow \hat P   \cr
                    &                                                        &               \cr 
                    &                                &      H^\C        \cr
\end{matrix}$$

\nmedskip
where $\, \hat \Omega_H \,$ is the envelope
of holomorphy of  $\, \Omega_H \,$. Moreover
$\, H\,$
acts on $\, \hat \Omega_H \,$ and all maps are 
$\, H$-equivariant. Furthermore 
$\,\hat P \,$ is injective by [CIT,\,Th.\,p.\,543] and if
$\, \Omega_H \,$
is holomorphically convex then 
$\, j \,$ is  surjective and consequently  
$\, \hat \Omega_H \,$ is biholomorphic to 
$\, P(\Omega_H) \,$, giving a contradiction.
Hence $\,\Omega_H\,$ is not holomorphically convex.

Notice that $\,\hat \Omega \cong \hat P(\hat \Omega)\,$ contains
$\, P(\Omega_H) = H\cdot \exp_{H^\C}(i\tilde P(O_0)) \,$ 
and the convex envelope
of $\, \tilde P(O_0) \,$ is $\,\h\,$, thus by the above arguments the envelope
of holomorphy $\,\hat \Omega \,$ is biholomorphic to $\,H^\C \cong
\C^3$.

Finally the maximal convex $\,\sigma$-invariant subdomain of 
$\, \tilde P(O_0)\,$ is $\,\{ \, (A,B,C) \in \h \, : \, A^2 +B^2 <
3\,\} = \tilde P(O_2)\,$, where $\,O_2:= \{ \, (a,b,c) \in \h \, : \,
a^2 +b^2 < 3\frac{c^2}{\sinh^2(c)}\,\}$. One checks that 
$\,O_2 \subset O_1\,$, thus $\,\tilde P|_{O_2}\,$ is injective and
$\, H \cdot O_2\,$ is biholomorphic to $\,H\cdot \exp_{H^\C}(i\tilde
P(O_2))\,$. Moreover one has 
\begin{align*}
\exp_{H^\C}(\,(A',B',C') + & i(A,B,C) \,)\,=  \\  
                & \exp_H(A',B',C')\,
\exp_{H^\C}(\,i(A,\,B,\,C-\frac{1}{2}(A'B-AB')\,)\,).
\end{align*} 

\noindent
It follows that
$$\,\exp_{H^\C}^{-1}(H\cdot \exp_{H^\C}(i\tilde P(O_2)) =
\,\{\, ({\mathrm{Im}}\,z_1)^2 + ({\mathrm{Im}}\,z_2)^2 < 3 \,\},$$

\nmedskip
where $\,(z_1, z_2,z_3)=(A'+iA,B'+iB,C'+iC)\,$ are natural complex 
coordinates of $\,h^\C \cong \C^3\,$ and this yields the statement. 
 \end{proof}

\bigskip
\begin{rem}
\label{O_1image}
Since $\,\tilde P \,$ is injective
on $\,O_1\,$, then 
the $\,H$-invariant domain
defined by $O_1$ is holomorphically separable.
As a matter of fact one may 
show that $\,\tilde P(O_1)  
\,=\, \h \setminus \{ \,(A,B,C) \in \h
\ :\  A^2 +B^2 \geq 3, \ C=0 \}$
and  analogous arguments as above
show that such $\,H$-invariant domain
is not holomorphically convex.
\end{rem}

\bbigskip


\section{Generalized Heisenberg groups}

\bigskip

Here we apply  results of the previous section
to generalized Heisenberg
groups exhibiting  additional examples
of non-Stein maximal domains of existence for the adapted
complex structure. We refer to \cite{BTV}
for the basic properties of generalized 
Heisenberg groups. 

Let $\,G\,$ be a generalized Heisenberg
group with Lie algebra $\,\g\,$, consider its abelian subalgebra
$\,\mathfrak{z}:= [\mathfrak{g},\mathfrak{g}]\,$
and the subspace $\,\mathfrak{v}\,$
orthogonal to $\,\mathfrak{z}\,$ with respect to
the $\,G$-invariant metric
$\,(\,\cdot \, ,\,\cdot \, ) \,$ of $\,G\,$.
Then for all 
$\,V +Y \in  \mathfrak{v} \oplus \mathfrak{z} = \g \cong T_eG\,$
the unique geodesic tangent
to $\,V +Y \,$ at $\,0\,$ can
be explicitely given  by 
$ \,t \mapsto \exp_G \circ \varphi_G(t, V+Y)\,$
for a certain  real-analytic map
$\,\varphi_G: \R \times T_eG \to \g \,$
(see [BTV,\,Th.\,p.\,31]).

A straightforward computation shows that $\,\varphi_G\,$
extends holomorphically on $\,\C\times T_eG\,$
to $\,\g^\C\,$ and 
analogous arguments as in the previous section 
imply that 
$\,\Omega_G := G \cdot O_G \subset TG\,$
is the maximal the adapted complexification,
where $\, O_G \,$ is the connected component of
$$\{ \, V+Y  \in T_eG \ :\ \det(D\tilde P_G)_{V+Y} \not=0 \,\}$$

\nmedskip 
containing $\,0\,$ and
$\, \tilde P_G: T_eG \to \g \,$ is given
by

$$\, V+Y  \mapsto  l(|Y|)V \ + \  (\,1 + |V|^2m(|Y|)\,)Y \,.$$

\nmedskip
Here $\,|\, \cdot \,|\,$ denotes the norm induced by
$\,(\,\cdot \, ,\,\cdot \, ) \,$ 
and the real-analytic
functions $\,l,m: \R \to \R\,$
are   defined by
$$l(t) := \frac{\sinh(t)}{t}, \quad \quad \quad  m(t) :=
 \frac{t -\sinh(t)\,\cosh(t)}{2t^3}.$$

\nsmallskip
Now for $\,V+Y \in \mathfrak{v} \oplus \mathfrak{z}=T_eG\,$
with $\,Y \not= 0\,$ and $\,  U+X \in T_{V+Y}T_eG \cong T_eG\,$
one has 

$$(D \tilde P_G)_{V + Y}(U + X)  \,=  \,\frac{\partial}{\partial t}
\,\tilde P_G( (V + tU) + (Y + t X) ) \vert_{t = 0} \, =
\,  l^{\prime}(\vert Y \vert)
\frac{(Y,X)}{\vert Y \vert} V \,+\, l( \vert Y \vert) U $$
$$ + \ \left(2 (V, U) m(\vert Y \vert) + \vert V \vert^2 
      m^{\prime}(\vert Y \vert) \frac{ (Y,X) }{\vert Y \vert} \right)Y
\ + \  (1 + \vert V \vert^2 m(\vert Y \vert) )X.$$

\nsmallskip
Note  that the equation is  written according to the
splitting $\, \mathfrak{v} \oplus \mathfrak{z}\,$
and  since  $\,l(\vert Y \vert)\,$ never vanishes,
the $\,\mathfrak{v}\,$-part vanishes
if and only if 

$$ U \ =  \ 
- \,\frac{(Y,X) \,l^{\prime}(\vert Y \vert)}{\vert Y \vert
\, l(\vert Y \vert)} \, V.$$

\nmedskip
It follows that the central $\mathfrak{z}$-part also vanishes
if and only if 
$$ ( 1 + \vert V \vert^2 m(\vert Y \vert) ) X \ =\ 
    \vert V \vert^2\,\frac{(X, Y) }{\vert Y \vert}  
    \left( 2 \frac{ l^{\prime}(\vert Y \vert)}{l(\vert Y \vert)}
 m(\vert Y \vert)
   \ - \ {m^{\prime}(\vert Y \vert)}  \right) Y. $$

\noindent
In particular $\,Y\,$ and $\,X\,$ have to be proportional. 
Since both sides are homogeneous of degree 1
in $\,X\,$, then $\,(D \varphi) \vert_{(V + Y)}\,$
is singular if and only if 
\begin{equation}
\label{lowrank}
1 + \vert V \vert^2 m(\vert Y \vert)\, =\,
 \vert V \vert^2 \,\vert Y \vert \,  
    \left(  2 \frac{ l^{\prime}(\vert Y \vert)}{l(\vert Y \vert)}
 m(\vert Y \vert)\   - \ {m^{\prime}(\vert Y \vert)} \right)  
\end{equation} 

\noindent
An analogous computation shows that $\,(D \tilde P_G)_{V+Y}\,$
has maximal rank if $\,Y=0\,$, thus equation (\ref{lowrank})
describes the singular locus of $\,D \tilde P_G \,$.
It is remarkable that this identity is independent of the
fine structure of the generalized Heisenberg group, e. g.  of
its dimension or the dimension of its centre.
In particular if $\,H\,$ is the 3-dimensional Heisenberg 
group considered in the previous
section, equation (\ref{lowrank}) determines
the border of $\,O_H = \Omega_H \cap T_eH\,$.
Using this fact we are now going to show
that for a generalized Heisenberg group
$\,G\,$ there exist many copies of $\,\Omega_H\,$
embedded as closed submanifolds in $\,\Omega_G\,$.

\medskip
Let $\,G\,$ be 
a generalized Heisenberg group and choose 
non zero elements
$\,\bar V_1 \in \mathfrak v\,$ and 
$\,\bar Y \in \mathfrak z\,$.
Then there exists an element 
$\,\bar V_2 \in \mathfrak v\,$ such that the closed subgroup
$\,\exp_G(\mathrm{span}\{\bar V_1,\bar V_2,\bar Y\})\,$ 
is a totally geodesically embedded 3-dimensional
Heisenberg group (see [BTV,\,p.\,30]). 
Denote by $\, I:H \to G\,$ such an embedding
and note that since
$\,\exp_{G^\C}: \g^\C \to G^\C \,$ is a biholomorphism,
then $\,I\,$ extends to a holomorphic embedding $\,I^\C:H^\C
\to G^\C\,$ of the universal complexification 
of $\,H\,$ into the universal complexification 
of $\,G\,$ such that the diagram
$$\begin{matrix}  
\h^\C \ \ \ \ \ \   &   \buildover{DI^\C}   \to  &   \g \ ^\C \ \ \ \  \cr
                                    &                              &           \cr
         \downarrow\, \exp_{H^\C}   &         &     \ \downarrow \, \exp_{G^\C}  \cr
                           &                               &                   \cr
         H ^\C\ \ \ \ \ \ \  & \buildover{I^\C}   \to  \    & G^\C \ \ \ \ \ \cr
\end{matrix} 
$$

\nmedskip
commutes.  Now $\,I:H \to G\,$ is totally geodesic, thus
$\  t \mapsto I \circ \exp_H \circ \varphi_H(t,v)\ $ is the 
unique geodesic of $\,G\,$ tangent to $\,DI(v)\,$ at $\,0\,$
for all $\,v \in \h \,$. Then 
$$I \circ exp_H \circ \varphi_H(t,v) \,=\,
exp_G \circ \varphi_G(t,DI(v))$$

\nsmallskip
and by the identity principle
\begin{equation}
\label{diagram2}
\,I^\C \circ \exp_{H^\C}\circ \varphi_H(z, v)
\, =\, \exp_{G^\C}\circ \varphi_G(z, DI(v))
\end{equation}

\nsmallskip 
for all $\,z\in \C\,$, since this holds for all $\,z\in \R$.
Commutativity of the diagram

\begin{equation}
\label{final}
\begin{matrix}  
TH  \ \ \ \   &   \buildover{DI}   \to  &   TG\ \ \ \  \cr
                                    &                              &           \cr
         \downarrow\, P_H   &          &     \ \downarrow \, P_G  \cr
                           &                               &                   \cr
       H^\C  \ \ \ \ \  & \buildover{I^\C}   \to  \    &  G^\C  \ \ \ \cr
\end{matrix} 
\end{equation}

\nmedskip 
follows. For this note that being $\,I\,$ a group
homomorphism then $\,DI:TH \to TG\,$ is
$\,H$-equivariant, i.e., $\,DI(g_*w) = I(g)_*DI(w)\,$
for all $\,g \in H\,$ and $\,w \in TH\,$. In particular
$$P_G \circ DI(g_*v) \,=\, P_G (\, I(g)_*DI(v)\,) =
I(g)\,\exp_{G^\C}\circ \varphi_G(i,DI(v))$$

\nmedskip
for all $\,g \in H \,$ and $\,v \in T_eH \,$.
On the other hand using equation (\ref{diagram2}) one
obtains
$$I^\C \circ P_H (g_*v) \,=\, I^\C (g\,
\exp_{H^\C} \circ \varphi_H(i , v)\,)\,=\,
\, I^\C (g)\,I^\C(\,
\exp_{H^\C} \circ \varphi_H(i , v)\,) \,=$$
$$I(g)\,\exp_{G^\C}\circ \varphi_G(i , DI(v)),$$

\nmedskip
showing that the above diagram is commutative.

From the equivariance of $\,DI\,$ it follows that
$\,DI(\Omega_H) =\,DI(H\cdot O_H)\,=\,I(H) \cdot DI(O_H)\,$
and since $\,DI\,$ is isometric and 
the border of $\,O_H\,$ is defined by
equation (\ref{lowrank}) which also
describes the singular locus of  $\,\tilde P_G\,$ one has 
$$DI(O_H) \,\subset \, O_G \cap 
\mathrm{span}\{\bar V_1,\bar V_2,\bar Y\},
\quad \quad \quad \quad DI(\partial O_H) \subset \partial O_G.$$

\nsmallskip
Thus $\,DI(\Omega_H)
\,\cong I(H) \times DI(O_H)\,$
is closed in $\,\Omega_G \cong G  \times O_G\,$. 

Furthermore $\,DI\,$ is injective and $\,P_H\,$, $\,P_G\,$
are locally biholomorphic
where the adapted complex structure
is defined (cf. Theorem\,\ref{sliceflow}), thus diagram (\ref{final}) shows that
$\,DI(\Omega_H) \cong \Omega_H\,$ is a closed
complex submanifold of $\,\Omega_G\,$. Finally
by Proposition 
 (\ref{notstein})   the domain $\,\Omega_H\,$
is neither holomorphically separable, nor holomorphically
convex, thus one has 

\bigskip
\begin{prop}
\label{generalized} Let $\,G\,$ be a generalized 
Heisenberg group. Then the maximal adapted
complexification
$\,\Omega_G\,$ is neither  holomorphically
separable, nor holomorphically convex.
\end{prop}

\bigskip
\medskip

\bbigskip

Stefan Halverscheid

Math. Fakult\"at, Ruhr-Universit\"at Bochum

Universit\"atsstr. 150, D-44780 Bochum, Germany

{\tiny {\it E- mail}: sth@cplx.ruhr-uni-bochum.de}

\medskip
\bigskip

Andrea Iannuzzi 

Dip. di Matematica, Universit{\`a} di Bologna

P.$\,$zza di Porta S. Donato 5, I-40126 Bologna, Italy

{\tiny {\it E-mail}: iannuzzi@dm.unibo.it}

\end{document}